\documentclass{tran-l}
\usepackage{amssymb,latexsym, amscd}
  \usepackage[all]{xy}
 \newtheorem{theorem}{Theorem}[subsection]
 
 \newtheorem{lemma}[theorem]{Lemma}
 
 \theoremstyle{definition}
 \newtheorem{definition}[theorem]{Definition}
 \theoremstyle{definition}
 
 \theoremstyle{remark}
 
 \numberwithin{equation}{subsection}

\usepackage{latexsym}
\usepackage{amssymb}

\newcommand{\ben}{\begin{equation}}
\newcommand{\een}{\end{equation}}

\newcommand{\integer}{\ensuremath{{\mathbb Z}}}

\newcommand{\complex}{\ensuremath{{\mathbb C}}}

\newcommand{\GL}[1]{\ensuremath{{\mathrm {GL}_{ #1 }}}}
\newcommand{\GLC}[1]{\GL{#1}(\complex)}
\newcommand{\SL}[1]{\ensuremath{{\mathrm {SL}_{ #1 }}}}
\newcommand{\SLC}[1]{\SL{#1}(\complex)}

\newcommand{\XX}{{\mathcal X}}

\newcommand{\VV}{{\mathcal V}}

\newcommand{\LL}{\mathcal{L}}

\newcommand{\Gg}{\mathsf{G}}
\newcommand{\Dd}{\mathsf{D}}

\newcommand{\Ll}{\mathsf{L}}

\newcommand{\Loop}{\mathsf{L}}

\newcommand{\timests}{\: {}_{t}  \! \times_{s}}

\newcommand{\morb}{\mu^{\mathrm{orb}}}
\newcommand{\Tate}{\ensuremath{{\mathbb L}}}
\newcommand{\Disc}{\ensuremath{{\mathbb D}}}
\newcommand{\Spec}{\ensuremath{{\mathrm{Spec}}}}
\newcommand{\Conj}{\ensuremath{{\mathrm{Conj}}}}
\newcommand{\rank}{\ensuremath{{\mathrm{rank}}}}

\begin{document}

\title[Global McKay-Ruan Correspondence via Motivic Integration]{The Global McKay-Ruan Correspondence via Motivic Integration}
\author{ Ernesto Lupercio and Mainak Poddar}

\address{Department of Mathematics, University of Wisconsin at Madison, Madison, WI 53706}
\address{Department of Mathematics, Michigan State University, East Lansing,
 MI 48824} \email{ lupercio@math.wisc.edu \\ poddar@math.msu.edu}


\begin{abstract}
The purpose of this paper is to show how the methods of motivic
integration of Kontsevich, Denef-Loeser
 and Looijenga
can be adapted to prove the McKay-Ruan correspondence, a
generalization of the McKay-Reid correspondence to orbifolds that
are not necessarily global quotients.
\end{abstract}

\maketitle

\section{Introduction}

\subsection{}
  In this paper when we say orbifold we mean an algebraic variety with quotient
singularities. We will always work over the field $\complex$. We
say that an orbifold $X$ is a \emph{global quotient} if it is of
the form $X=U/G$ for $U$ smooth and $G$ a finite group. There are
simple examples of orbifolds that are not global quotients such as
weighted projective spaces like, for example, $WP(1,1,2)$. An
$n$-dimensional orbifold is \emph{Gorenstein} or
 $SL$ if all the local isotropy groups are finite subgroups of $\SLC{n}$.
 Given two Gorenstein or
 $SL$ orbifolds $X$ and $Y$ we say that they are
\emph{$K$-equivalent} if there is a common birational resolution
$\phi\colon Z \to X$, $\psi \colon Z \to Y$ so that $\phi^* K_X =
\psi^* K_Y$.

\subsection{}
 The cohomological \emph{McKay-Reid} correspondence
 states that for an $SL$ global quotient orbifold having
 a smooth crepant resolution, there is a certain correspondence
 between the cohomology generators of the resolution and
 the conjugacy classes of the group $G$. A very good introduction  to this subject, including its relation to the classical McKay correspondence can be found in \cite{Reid}.  We refer the reader to \cite{Batyrev2,
 DenefLoeser2, ItoReid, Kaledin}
 for detailed results. Yongbin Ruan (cf. \cite{RuanString} Section 6.4)
 conjectured a generalization
 of this correspondence to general orbifolds, which has the added
 advantage that the smoothness of the resolution is no longer essential.
 In the \emph{McKay-Ruan} correspondence (stated below as Theorem \ref{McKay}), the role of the conjugacy
 classes is taken over by the \emph{twisted sectors}  of the orbifold in question (see \cite{ChenRuan,
 Dixon, Kawasaki} and Definition \ref{inertia} below).

\subsection{} The purpose of this note is to
 show how the methods of motivic
integration of Kontsevich \cite{Kontsevich}, Denef-Loeser
\cite{DenefLoeser1, DenefLoeser2} and Looijenga \cite{Looijenga}
can be adapted to prove the \emph{McKay-Ruan} correspondence,
namely the following theorem.

\begin{theorem}\label{MRuan}
If the orbifolds $X$ and $Y$ are $K$-equivalent and complete then their
orbifold Hodge numbers, orbifold Hodge structures and orbifold
Euler characteristics coincide.
\end{theorem}

This theorem has been proved independently by T. Yasuda
\cite{Yasuda} as we were informed by an email from him during the
preparation of this note.

\subsection{} The previous theorem is actually a consequence of the
transformation rule for the motivic measure (Theorem 1.16 in
\cite{DenefLoeser2}, Theorem 3.2 in \cite{Looijenga}), and the
following formula that is analogous to Theorem 8.1 in
\cite{Looijenga} except that in there the orbifold was only a global quotient.

\begin{theorem}\label{McKay}
$$ \morb(\LL(X))  = \sum_{\alpha \in \pi_0(\wedge X)} [X^\alpha /
 X]\Tate^{w(\alpha)} \in \hat{M}_X[\Tate ^{1/m}] $$
\end{theorem}

Let us briefly describe the different terms involved in this
formula (cf.~\cite{Looijenga}). $\LL(X)$ is the space of arcs in
$X$ whose $\complex$-points correspond to formal arcs $\Disc \to
X$. Here $\Disc = \Spec \complex[[z]]$. The \emph{motivic ring}
$M_X$ is equal to the localization of $K_0(\VV_X)$, the
Grothendieck ring of $X$-varieties (plus the relation 1.10 in
\cite{DenefLoeser2}), with respect to the Tate motive
$\Tate=[\mathbb{A}^1]\in K_0(\VV_X)$. The $M_X[\Tate ^{1/m}]$
valued measure $\morb$ on $\LL(X)$ is defined by Denef-Loeser in
\cite{DenefLoeser2} \S 2.7 (cf. \cite{Looijenga} \S8). This
measure is determined by the dualizing sheaf of $X$ if $X$ is
Gorenstein. We write $\wedge X$ to denote the twisted sectors or
\emph{inertia orbifold} of $X$ (see Definition \ref{inertia}, cf.
\cite{LupercioUribe} definition 3.6.5.). Finally the numbers
$w(\alpha)$ are the \emph{degree shifting numbers} described by
Chen-Ruan \cite{ChenRuan}. They are also known as \emph{Fermionic
degree shifting numbers} in the  physics literature and sometimes
referred to as \emph{age} in the terminology of Miles Reid.

\section{The Case of the Global Quotient.}

\subsection{}
 The proof of Theorem \ref{McKay} in the global quotient case can
be found in \cite{DenefLoeser2, Looijenga}. The proof in the
general case follows closely the proof of the McKay-Reid
correspondence (Theorem 8.1 in \cite{Looijenga}, cf. \cite{Batyrev1,
DenefLoeser2}) described by Looijenga \cite{Looijenga}, Section 8.
We briefly recall the basic idea of this argument now (cf.
\cite{Reid} formula 4.5). We divide the proof into two steps.

\subsection{Step I}\label{StepI} Let $p \colon U \to X$ be the quotient map associated to
$X$. We set $\LL^\circ(X)$ to be the set of arcs in $X$ not
contained in the discriminant of $p$ (this is written $\LL'(X)$ in
\cite{Looijenga} and $\LL^g(X)$ in \cite{DenefLoeser2}).

We denote by $[m]\colon\Disc \to \Disc$ the $m$-th power map
$z^{1/m} \mapsto z$. Let $\zeta_m=\exp(2\pi i/m)$ be a primitive
$m$-root of unity where $m$ is the order of $G$.

 The first thing to do is to verify the
following decomposition (cf. \cite{Reid} (4.6) ,
\cite{DenefLoeser2} 6.2 and \cite{LupercioUribe} \S 6.1.2)
$$ \LL^\circ(X) = \coprod_{(g)\in\Conj(G)} \LL^\circ_g (X). $$
Here $\LL^\circ_g(X)$ consists of the arcs $\gamma$ in
$\LL^\circ(X)$ that have a  lift $\tilde{\gamma}$ in $\LL^\circ(U)$
with the following property
$$ g \tilde\gamma = \gamma \zeta_m .$$
We call $\LL^\circ_g(X)$ the arcs on $X$ twisted by $g$.

\subsection{Step II} The next step is to verify the following identity,
$$ \morb(\LL_g (X))  =  [X^g / X]\Tate^{w(g^{-1})}  $$
where $X^g = U^g / C(g)$ is the $(g)$-twisted sector, namely the
fixed points of $g$ modulo its centralizer.

This calculation is performed using the so-called change of
variables formula (see \cite{DenefLoeser2} Theorem 1.16,
\cite{Craw} Theorem 2.18, \cite{Looijenga} Theorem 3.2) for the
motivic integral defining the orbifold motivic measure. This is
done for example in \cite{DenefLoeser2} \S 3.

\subsection{} Finally all there is to do is to add over all
$(g)\in\Conj(G)$ to get
$$ \morb(\LL(X))= \morb(\LL^\circ (X)) = \sum_{(g)} [X^g/X]\Tate^{w(g^{-1})}. $$
Since $\pi_0 ( \wedge[U/G] ) = \Conj(G)$ (see
\cite{LupercioUribe} Proposition 6.2.1) we have Theorem
\ref{McKay} for the case $X=U/G$.

\section{The General Case.}

\subsection{} We can think of a general orbifold $X$ as a Deligne-Mumford
stack and such a stack admits an atlas given by an \emph{\'{e}tale
separated groupoid $\Gg$ in schemes.} The category of orbifolds is
then equivalent to the category of groupoids up to Morita
equivalence. We refer the reader to \cite{Moerdijk, Pronk,
Abramovich, LupercioUribe1, LupercioUribe, Totaro} for details. In
any case we will denote by $U=\Gg_0$ and $R=\Gg_1$ the smooth
schemes of objects and morphisms (arrows) of the groupoid $\Gg$
respectively, and the structure maps by:
       $$\xymatrix{
         \Gg_1 \timests \Gg_1 \ar[r]^{m} & \Gg_1 \ar[r]^i &
         \Gg_1 \ar@<.5ex>[r]^s \ar@<-.5ex>[r]_t & \Gg_0 \ar[r]^e & \Gg_1
         }$$
where $s$ and $t$ are the source and the target maps of morphisms,
$m$ is the composition of two of them whenever the target of the
first equals the source of the second, $i$ gives us the inverse
morphism and $e$ assigns the identity arrow to every object. We
will write $\XX=[U/\Gg]$ for the associated stack and $X=U/\Gg$
for the corresponding coarse moduli space. As before we write $p
\colon U \to X $ to denote the quotient morphism.

\subsection{The Arc Groupoid.}\label{arcgroupoid} The following definitions are
borrowed from \cite{LupercioUribe}. Let $\Disc_k = \Spec
\complex[z]/z^{k+1} $. We fix once and for all a positive integer
$m$ so that the orders of all the stabilizers of $X$ divide $m$
(here we are assuming that $X$ is complete). As before consider
the $m$-th power morphism $[m]:\Disc_k \to \Disc_k$ given by
$z^{1/m} \mapsto z$ and $\zeta_m = \exp(2\pi i/m)$. We denote by
$C_m$ the cyclic group generated by $\zeta_m$. Let $\Dd_k$ be the
groupoid with $(\Dd_k)_0 = \Disc_k$, $(\Dd_k)_1 = \Disc_k\times
C_m$, $s=id$ and $t$ being the Galois action of $\zeta_m$ in
$\Disc_k$.

\subsubsection{} We define a \emph{$k$-jet on the groupoid $\Gg$} to
be a morphism of groupoids
$$\gamma \colon \Dd_k \rightarrow \Gg.$$
We write $\gamma_i \colon (\Dd_k)_i \rightarrow \Gg_i$ for $i=0,1$
for the corresponding morphisms in objects and arrows of $\Dd_k$.

\subsubsection{} We define the \emph{$k$-jet groupoid} $\Loop_k \Gg$
associated to $\Gg$  by the following data:
\begin{itemize}
\item[(i)] Objects ($(\Loop_k \Gg)_0$):
 Morphisms $\Dd_k \to \Gg$.
\item[(ii)] Morphisms ($(\Loop_k \Gg)_1$): For two elements in $(\Loop_k \Gg)_0$, say
$\Psi, \Phi : \Dd_k \to \Gg$  a morphism (arrow) from $\Psi$ to
$\Phi$ is a morphism $\Lambda : \Disc_k \times C_m \to \Gg_1$ that
makes the following diagram commute
      $$
      \xymatrix{
       \Disc_k \times C_m \ar[r]^\Lambda \ar[d]_{s\times t}
      & \Gg_1 \ar[d]^{s \times t}\\
       \Disc_k \times \Disc_k \ar[r]_{(\Psi_0 , \Phi_0)} & \Gg_0 \times \Gg_0 }
      $$
and such that for $r \in \Disc_k \times C_m $
$$\Lambda(r)=\Psi_1(r) \cdot \Lambda(es(r))= \Lambda(et(r)) \cdot \Phi_1(r).$$
\end{itemize}

The composition of morphisms is defined pointwise, in other
words, for $\Lambda$ and $\Omega$ with
      $$\xymatrix{
      \Psi \ar@/^/[r]^\Lambda  & \Phi \ar@/^/[r]^\Omega & \Gamma}
      $$
we set $$\Omega \circ \Lambda(es(r)) := \Lambda(es(r)) \cdot
\Omega(es(r))$$ and
 $$\Omega \circ \Lambda(r):= \Omega \circ \Lambda (es(r)) \cdot \Gamma(r)=\Psi(r) \cdot
\Omega \circ \Lambda (et(r)).$$

\subsubsection{} The scheme structure on the space of objects
$(\Loop_k \Gg)_0$ is given by identifying it as a subscheme of
$\LL_k(\Gg_1)$, for $\gamma_1$ completely determines $\gamma$.
Similarly $(\Loop_k \Gg)_1$ is naturally a subscheme of
$\LL_k(\Gg_1)^m$. In fact we can do better. By 3.2.4 in
\cite{LupercioUribe} or simply by recalling that $\Gg$ is
\'{e}tale  we can show as in 3.2.6\cite{LupercioUribe} that
$\Loop_k \Gg$ is actually an \'{e}tale groupoid in schemes, but we
will not need this.

\subsubsection{} The \emph{arc groupoid} is similarly defined. We
denote it by $\Loop \Gg = \Loop_\infty\Gg$. Just as in
\cite{Looijenga} for $m \geq n$ we have projections $\pi_n^m
\colon \Loop_m\Gg \to \Loop_n\Gg$. We will simply write
$\pi_n$ for $\pi_n^\infty$.

\subsection{Step I} Now we generalize \ref{StepI} to the general
case.

\begin{definition}\label{inertia} The inertia groupoid $\wedge \Gg$ is defined
in the following way:
\begin{itemize}
\item[(i)] Objects $(\wedge \Gg)_0$: Elements $v \in \Gg_1$ such that $s(v) = t(v)$.
\item[(ii)] Morphisms $(\wedge \Gg)_1$: For  $v,w \in (\wedge \Gg)_0$ an
arrow $v \stackrel{\alpha}{\to} w$ is an element $\alpha \in
\Gg_1$ such that $v \cdot \alpha = \alpha \cdot w$
       $$
        \xymatrix{
        \circ \ar@(ul,dl)[]|{v} \ar@/^/[rr]|{\alpha}
        &&\circ \ar@(dr,ur)[]|{w^{-1}} \ar@/^/[ll]|{\alpha^{-1}}
        }$$
\end{itemize}
\end{definition}

It is known that the inertia groupoid in the case of an orbifold
matches with what is commonly known in the literature as twisted
sectors \cite{LupercioUribe}.

The inertia groupoid defines a stack because as Moerdijk points
out \cite{Moerdijk} it can be seen as $$\wedge \Gg = S_\Gg
\rtimes \Gg$$ where $S_\Gg = \Delta^*(\Gg_1)$ and
$\Delta\colon \Gg_0 \to \Gg_0 \times \Gg_0$ is the diagonal
morphism.

\subsubsection{} Notice that it is clear from the definitions that
$\Loop_0(\Gg)= Hom(\overline{\integer}, \Gg) = \wedge\Gg$ (cf.
3.6.4\cite{LupercioUribe}.) We will also call the map $\pi_0$ the
\emph{evaluation map} $ev\colon \Ll \Gg \to \wedge \Gg$. We will
write $\tau$ to denote the composition
$$\tau\colon\Ll \Gg \rightarrow\wedge \Gg \rightarrow\pi_0(\wedge \Gg )$$
and will write for $\alpha \in \pi_0 (\wedge \Gg)$
$$\Ll_\alpha(\Gg) = \tau^{-1} (\alpha)$$

\subsubsection{} It may be worth pointing out here for the sake of comparison with the
global quotient case that when $X=U/G$ we have (6.2.2
\cite{LupercioUribe})
$$\Conj(G) = \pi_0 (\wedge(X))$$

\subsubsection{} Let $p \colon U \to X$ be the quotient map associated to
$X$. We set $\LL^\circ(X)$ to be the set of arcs in $X$ not
contained in the discriminant of $p$. (This is written $\LL'(X)$ in
\cite{Looijenga} and $\LL^g(X)$ in \cite{DenefLoeser2}.) Just as
before the measure of $\LL(X)-\LL^\circ(X)$ is zero. (This is a
local statement after all, but in any case it follows just as in
2.1 \cite{DenefLoeser2}.) Compare with section 8 in
\cite{Looijenga}.

\subsubsection{} We write $\LL_\alpha (U)$ to denote the arcs in $U$ of the form
$\gamma_1$ for $\gamma$ an object in the jet groupoid lying in
$\tau^{-1} (\alpha)$ (in the global quotient case this is what we
called $\LL_g (U)$.) Given an arc in $\LL^\circ(X)$ the map
$\gamma\circ[m]$ (where $[m]$ is defined in \ref{arcgroupoid})
lifts to a morphism $\Dd \to \Gg$. Let $\LL^\circ_\alpha(X)$ the
set of arcs corresponding to $\Loop_\alpha \Gg$ under this lift
(cf. \cite{Looijenga} discussion after 8.3).

\subsubsection{} We have the following decomposition $$
\LL^\circ(X) = \coprod_{\alpha \in \pi_0 (\wedge \Gg)} \LL^\circ_\alpha (X). $$ This
is true by \S 2.1 \cite{DenefLoeser2}, for again this is a local
statement.

\subsection{Step II} We want to compute $\morb(\LL^\circ_\alpha (X))$.

\subsubsection{} We will need the following stratification of an orbifold
due to Haefliger \cite{Haefliger} Proposition A.2.2.

\begin{lemma} For an orbifold $X$ that is complete, there exist a
stratification
$X=\coprod_{i=1}^N X_i$ and a corresponding stratification of the \'{e}tale
groupoid $\Gg=\coprod_i \Gg_i$ so that each $X_i$ is smooth with
stabilizer $\Gamma_i$ (where $\Gamma_i$ is a finite group) and
each $\Gg_i$ is Morita equivalent to the action groupoid $X_i
\times \Gamma_i \rightrightarrows X_i$.
\end{lemma}

\begin{proof} By the results of \cite{Totaro} and
\cite{Luna} it is enough to consider the case when $X=M/\GLC{n}$
is given by the action groupoid $M\times \GLC{n}
\rightrightarrows M$. Make a list of all possible stabilizers
$G_i^k$ for $i=1,\ldots,N$, so that for a fixed $i_0$ all the
$G_{i_0}^k$ are conjugate in $\GLC{n}$. Let
$N_i^k=N_{\GLC{n}}(G_i^k)$ be the normalizer of $G_i^k$,
$H_i^k=N_i^k/G_i^k$, and set $\Gamma_i = G_i^{k_i}$ for a fixed
$k_i$. Let $M_i^k$ the subset of $M$ with stabilizer equal to
$G_i^k$ (which is clearly constructible) and $M_i=\coprod_k
M_i^k$. The groupoid $M_i\times G \rightrightarrows M_i$ is a
subgroupoid of $M\times G \rightrightarrows M$. Note that $M_i\times
G \rightrightarrows M_i$ defines the same orbifold as $M_i^{k_i}
\times N_i^{k_i} \rightrightarrows M_i^{k_i}$. Set
$X_i=M_i^{k_i}/H_i^{k_i}$ (which is smooth because $H_i^k$ acts
freely on $M_i^k$), then the groupoid $M_i^{k_i} \times N_i^{k_i}
\rightrightarrows M_i^{k_i}$ defines the same orbifold as $X_i
\times \Gamma_i \rightrightarrows X_i$.
\end{proof}

\subsubsection{} We will restrict our attention to a fixed stratum
$X_i$ of $X$ and consider the arcs based at $X_i$ that we will
write $\LL(X)_{X_i}$ (cf. \S2.1\cite{DenefLoeser2}). We will
follow Denef-Loeser \S3 \cite{DenefLoeser2} to compute
$\morb(\LL^\circ_{\alpha}(X)_{X_i})$, and obtain the desired result
by summing over strata.

\subsubsection{} Write
$\Gg_i = [R_i \rightrightarrows U_i]\simeq [X_i \times \Gamma_i \rightrightarrows X_i]$ to
denote the \'{e}tale groupoid corresponding to $X_i$ and $p_i\colon U_i \to
X_i$ the corresponding quotient map. Let $\nu_i$ denote the
normal bundle of $U_i$ in $U$.

\subsubsection{} Since $\LL^\circ_\alpha(X)_{X_i}$ depends only on a
formal neighborhood $V_i$ of $X_i$ in $X$ and by decomposing the
normal bundle $\nu_{i} = \oplus \nu_{i,\alpha}^k$ (cf. Lemma 8.4
\cite{Looijenga}, Lemma 2.2 \cite{DenefLoeser2}) on eigenspaces
for a generator $g_\alpha$ of $\alpha$ in $\Gamma_i$ (the group
$\Gamma_i$ acts trivially on $U_i$ and therefore linearly by
fibers on $\nu_i$) we can define $\tilde{\lambda} \colon
\nu_{\complex[z]} \to X \otimes \complex[z]$ by the formula
2.3\cite{DenefLoeser2} used fiberwise. We define $w(\alpha) =
\sum_k (1-k/m) \rank(\nu_{i,\alpha}^k)$ (as in
8.4\cite{Looijenga}). This number is independent of $i$ (cf.
\cite{ChenRuan}). By 2.3\cite{DenefLoeser2} we have that there is
a $\complex[z]$-morphism $\tilde{\lambda}_* \colon \LL(\nu_i) \to
\LL(X)_{X_i}$ with $\LL^\circ_\alpha(X)_{X_i} =
\tilde{\lambda}_*(\LL(\nu_i))\cap\LL^\circ(X)$ showing that
$\LL^\circ_\alpha(X)_{X_i}$ is a $\complex[z]$-semi-algebraic
subset of $\LL(X)$. Moreover we also have a $\complex[z]$-morphism
$\lambda_* \colon \LL(\nu_i/\Gamma_i) \to \LL(X)$ inducing a
bijection  $\lambda_*\colon \LL(\nu_i)/\Gamma_i \simeq
(\LL(\nu_i)/\Gamma_i) \cap \lambda_*^{-1}(\LL^\circ(X))$.
Proceeding as in 3.3\cite{DenefLoeser2} we get
$\morb(\LL^\circ_\alpha(X)_{X_i}) = \Tate^{w(\alpha)}
\morb_{\LL(\nu_i/R_i)}(\LL(\nu_i)/R_i)$. But we have as in
3.4\cite{DenefLoeser2} that
\begin{eqnarray*}
\morb_{\LL(\nu_i /R_i)}(\LL(\nu_i)/R_i)
& = & [\pi_r(\LL(\nu_i))/R_i]\Tate^{-(r+1)\rank(\nu_i)/m} + R'_M \\
& = &
[X^\alpha_i]\Tate^{(r+1)\rank(\nu_i)/m}\Tate^{-(r+1)\rank(\nu_i)/m}
+ R'_M
\end{eqnarray*}
where $\lim_{M\to\infty} R'_M = 0$. This proves that
$$ \morb(\LL^\circ_{\alpha}(X)_{X_i}) = [X^\alpha_i]\Tate^{w(\alpha)},$$
finishing the proof of theorem \ref{McKay}.

\subsubsection{} Finally by applying $\chi_h$, $\chi_{hn}$ and
$\chi_{top}$  (cf. \S2,\S4\cite{Looijenga}) to
$\morb(\LL(X))=\morb(\LL(Y))$
we obtain theorem \ref{MRuan}.

\subsection{} We would like to thank very enlightening
conversations with R. Kulkarni,
 T. Nevins, M. Reid, Y. Ruan, G. Segal and B. Uribe
during the preparation of this work. The second author would also
like to thank M. Artin and D.A. Cox for useful email
correspondences.

\bibliographystyle{amsplain}
\bibliography{motivicmckay}

\end{document}